\documentclass[12pt,a4paper]{amsart}

\usepackage{amsfonts,amsthm,amsmath,amssymb,amscd,mathrsfs}
\usepackage{amsaddr}
\usepackage{graphics}
\usepackage{indentfirst}
\usepackage{cite}
\usepackage{bm, enumerate,xcolor}
\usepackage[dvips]{epsfig}
\usepackage{latexsym}
\usepackage{float}

\newtheorem{theorem}{Theorem}[section]
\newtheorem{remark}{Remark}[section]
\newtheorem{definition}{Definition}[section]
\newtheorem{lemma}[theorem]{Lemma}
\newtheorem{corollary}[theorem]{Corollary}
\newtheorem{pro}[theorem]{Proposition}

\def\be{\begin{eqnarray}}
	\def\ee{\end{eqnarray}}
\def\ba{\begin{aligned}}
	\def\ea{\end{aligned}}
\def\bay{\begin{array}}
	\def\eay{\end{array}}
\def\bca{\begin{cases}}
	\def\eca{\end{cases}}

\def\bt{\begin{theorem}}
	\def\et{\end{theorem}}
\def\bc{\begin{corollary}}
	\def\ec{\end{corollary}}
\def\bl{\begin{lemma}}
	\def\el{\end{lemma}}
\def\bp{\begin{proposition}}
	\def\ep{\end{proposition}}
\def\br{\begin{remark}}
	\def\er{\end{remark}}
\def\bd{\begin{definition}}
	\def\ed{\end{definition}}
\def\bpf{\begin{proof}}
	\def\epf{\end{proof}}
\def\bex{\begin{example}}
	\def\eex{\end{example}}
\def\bq{\begin{question}}
	\def\eq{\end{question}}
\def\bas{\begin{assumption}}
	\def\eas{\end{assumption}}
\def\ber{\begin{exercise}}
	\def\eer{\end{exercise}}

\def\XXint#1#2#3{{\setbox0=\hbox{$#1{#2#3}{\int}$ }
		\vcenter{\hbox{$#2#3$ }}\kern-.6\wd0}}

\begin{document}
\title[Maximum principles for Laplacian and fractional Laplacian]
{Maximum principles for Laplacian and fractional Laplacian with critical integrability}
\author{Congming Li${^1}$ and Yingshu L\"{u}${^2}$}
\thanks{${^1}$School of Mathematical Sciences, CMA-Shanghai, Shanghai Jiao Tong University, Shanghai 200240, China. Email: congming.li@sjtu.edu.cn.}
	\thanks{${^2}$ Corresponding author. School of Mathematical Sciences, Institute of Natural Sciences, Shanghai Jiao Tong University, Shanghai 200240, China. Email: yingshulv@sjtu.edu.cn.}
\date{}
\maketitle

\begin{abstract}

In this paper, we study maximum principles for Laplacian and fractional Laplacian with critical integrability. We first consider the critical cases for Laplacian with zero order term and first order term. It is well known that for the Laplacian with zero order term $-\Delta +c(x)$ in $B_1$, $c(x)\in L^p(B_1)$($B_1\subset \mathbf{R}^n$), the critical case for the maximum principle is $p=\frac{n}{2}$. We show that the critical condition $c(x)\in {L^{\frac{n}{2}}(B_1)}$ is not enough to guarantee the strong maximum principle. For the Laplacian with first order term $-\Delta +\vec{b}(x)$($\vec{b}(x)\in L^p(B_1)$), the critical case is $p=n$. In this case, we establish the maximum principle and strong maximum principle for Laplacian with first order term.

We also extend some of the maximum principles above to the fractional Laplacian. We replace the classical lower semi-continuous condition on solutions for the fractional Laplacian with some integrability condition. Then we establish a series of maximum principles for fractional Laplacian under some integrability condition on the coefficients. These conditions are weaker than the previous regularity conditions. The weakened conditions on the coefficients and the non-locality of the fractional Laplacian bring in some new difficulties. Some new techniques are developed.

\end{abstract}

\small{{\bf Key words:} Maximum principles; Laplacian; fractional Laplacian; critical integrability

}

\section{Introduction}
\quad ~ Maximum principles are fundamental tools in the study of partial differential equations. The classical maximum principle for harmonic functions and subharmonic functions can be traced to the work of Gauss \cite{Gauss}. Hopf established the classical strong maximum principle in \cite{Hopf1,Hopf2} which is a basic building block for the analysis of the second order elliptic partial differential equations. Later, various versions of maximum principles have been discussed by many researchers, as in Littman \cite{Littman}, Trudinger \cite{Trudinger1,Trudinger2,Trudinger3}, Protter and Weinberger  \cite{M.W}, Gidas, Ni and Nirenberg \cite{G-N-N}, Brezis and Lions \cite{BL}, Berestycki and Nirenberg \cite{B-N}, Berestycki, Nirenberg and Varadhan \cite{BNV}, Brezis and Ponce \cite{BP}, Pucci and Serrin \cite{p-s2,p-s}, Vitolo \cite{V}, and Cavaliere and Transirico \cite{c-t}.

In the last few decades, the Schr\"odinger operator $-\Delta +c(x)$ has attracted a lot of attention from many scientists, where $c(x)$ is a given potential in an open connected set $\Omega\subset\mathbf{R}^n$. The classical strong maximum principle states that under the certain condition on $c(x)$, if $u$ satisfies $-\Delta u+cu\geq 0$, $u\geq 0$, and $u(x_0)=0$ for some point $x_0\in \Omega$, then $u\equiv 0$ in $\Omega$. It is known that the strong maximum principle holds if $c(x)\in L^p(\Omega)$ for some $p>\frac{n}{2}$ (See Serrin \cite{Se}, Stampacchia \cite{Stam}, Trudinger \cite{Trudinger2}). Later, the assumption for $c(x)$ is weakened if adding some vanishing condition on $u$.  Ancona \cite{AA}, Brezis and Ponce \cite{BP} showed that $u\equiv 0\ a.e.$ in $\Omega$ if $c(x)\in L^1(\Omega)$, $c\geq 0 \ a.e.$ in $\Omega$, and there exists a quasi-continuous function related to $u$ vanishes on a set of positive $H^1$-capacity in $\Omega$, respectively. In \cite{OP}, Orsina and Ponce proved that for $p>1$ and $c(x)\in L^p(\Omega)$, $u(x)=0 \ a.e.$ in $\Omega$ if $u(x)$ satisfies vanishing condition $\lim_{r\rightarrow 0}\frac{1}{B(x,r)}\int_{B(x,r)}u(x)dx=0$ for every point $x$ in a compact subset of $\Omega$ with positive $W^{2,p}$ capacity. In \cite{BST}, Bertsch, Smarrazzo and Tesei gave a necessary and sufficient condition for the validity of the strong maximum principle in one dimension.

In this paper, we focus on whether the critical integrability condition for $c(x)$ can ensure the classical strong maximum principle for the Schr\"{o}dinger operator. The strong maximum principle for the Schr\"{o}dinger operator asserted in \cite{Se,Stam,Trudinger2} requires that $c(x)\in L^p(\Omega)$ for some $p>\frac{n}{2}$. Then what happens in the case $p\leq \frac{n}{2}$? It is well-known that the strong maximum principle fails for $p<\frac{n}{2}$. For instance, if $B_1\subset\mathbf{R}^n$, the function $u(x)=|x|^2$ satisfies $-\Delta u+cu=0$ in $B_1$ with $c(x)=\frac{2n}{|x|^2}$.

\emph{A natural question is that can one obtain the strong maximum principle for the Schr\"{o}dinger operator when $p=\frac{n}{2}$?}

One of the key results of this paper is to answer this question. We show that $c(x)\in L^{\frac{n}{2}}(B_1)$ is not enough to ensure the strong maximum principle for the Schr\"{o}dinger operator no matter how small $\|c\|_{L^{\frac{n}{2}}(B_1)}$ is.

More precisely, our first main result can be stated as follows.

\begin{theorem}\label{thm2}
Let $n\geq 3$. There exist $c_{\epsilon}(x)\in L^{\frac{n}{2}}(B_1)$ and  $u_{\epsilon}(x)\in H^1(B_1)\cap C(\bar{B}_1)$, $\epsilon>0$ such that
\begin{eqnarray}\label{2}
  -\Delta u_{\epsilon}(x)+c_{\epsilon}(x)u_{\epsilon}(x)= 0 \ \ \mbox{in} \ B_1,
 \end{eqnarray}
and\\
~~~~(i) $u_{\epsilon}(x)\geq 1$ on $\partial B_1$,\\
\ \ (ii) $u_{\epsilon}(0)=0$,\\
\ \  (iii)$\lim\limits_{\epsilon \rightarrow 0^+}\|c_{\epsilon}\|_{L^{\frac{n}{2}}(B_1)}=0$.
\end{theorem}

For completeness, we present the maximum principle for the Schr\"{o}dinger operator with critical case $p=\frac{n}{2}$ in $B_1$.

\begin{theorem}\label{thm1}
  Assume that $c(x)\in L^{\frac{n}{2}}(B_1) $, and $u(x)\in H^1(B_1)$($n \geq 3$)  is a weak solution of
   \begin{eqnarray}\label{1}
 \left\{\begin{array}{l} -\Delta u(x)+c(x)u(x)\geq 0 \ \ \ \mbox{in} \ B_1  \\
u(x) \geq 0 \ \ \ \mbox{on}\ \partial B_1.
\end{array}
\right.
 \end{eqnarray}
 There exists a positive constant $k(n)$ such that if
 $\|c^-\|_{L^{\frac{n}{2}}(B_1)}\leq k(n)$, then $u(x) \geq 0$ in $B_1$. Here $c^{-}(x)=-\mbox{min}\{c(x),0\}$.
 \end{theorem}

\begin{remark}\label{rem}
Note that the maximum principle for the Schr\"{o}dinger operator(Theorem \ref{thm1}) is not true if $n=2$. For instance, one can take $u_{\epsilon}(x)=-ln(\epsilon \left|ln|x|\right|)$ changing sign in $B_1$, $\epsilon>0$, with corresponding $\|c_{\epsilon}\|_{L^1(B_1)}\rightarrow 0$ as $\epsilon \rightarrow 0^{+}$.
\end{remark}

\begin{remark}
Actually, it follows from Remark \ref{rem} that the strong maximum principle for the Schr\"{o}dinger operator with critical integrability $c(x)\in L^{\frac{n}{2}}(B_1)$ is not true if $n=2$.
\end{remark}

Next, we give a refined version of the strong maximum principle for the Schr\"{o}dinger operator which is useful for analysis of PDEs in practice. 
Various refined versions of the strong maximum principle have been studied, see \cite{L.W.X,L.L.W.X} and references therein. Note that in \cite{L.W.X}, they required that the corresponding coefficient $c(x)$ is bounded.

The following theorem is not really new. For convenience and completeness, we present it here.

\begin{theorem}\label{thm3}
 Assume that $c(x)\in L^p(B_1)$($p>\frac{n}{2}$), and $u(x)\in H^1(B_1)$($n \geq 3$) is a weak solution of
   \begin{eqnarray}\label{3}
 \left\{\begin{array}{l} -\Delta u(x)+c(x)u(x)\geq 0 \ \ \ \mbox{in} \ B_1  \\
u(x) \geq m>0 \ \ \ \mbox{on}\ \partial B_1.
\end{array}
\right.
 \end{eqnarray}
 There exists a positive constant $k(n,p)$ such that if
 $\|c^+\|_{L^{p}(B_1)}\leq k(n,p)$, then
  \begin{eqnarray*}
  u(x) \geq \gamma m \ \ \mbox{in} \ B_1,
 \end{eqnarray*}
 where $c^{+}(x)=\mbox{max}\{c(x),0\}$, and $\gamma$ is a positive constant depending only on $n$, $p$ and $\|c^+\|_{L^p(B_1)}$.
 \end{theorem}

In the following, we consider the maximum principle and strong maximum principle for Laplacian with first order term. Actually, the critical integrability condition for $\vec{b}(x)$ is $\vec{b}\in L^p(B_1)$ for $p=n$. The crucial observation is that compared to Theorem \ref{thm3}, the strong maximum principle for Laplacian with first order term is true in the critical case $p=n$.

\begin{theorem}\label{thm4}
   Assume that $\vec{b}(x)\in L^n(B_1)$, and $u(x)\in H^1(B_1)$($n \geq 3$) is a weak solution of
   \begin{eqnarray}\label{4}
 \left\{\begin{array}{l} -\Delta u(x)+\vec{b}(x)\cdot \nabla u(x)\geq 0 \ \ \ \mbox{in} \ B_1  \\
u(x) \geq 0 \ \ \ \mbox{on}\ \partial B_1.
\end{array}
\right.
 \end{eqnarray}
 There exists a positive constant $k(n)$ such that if
 $\|\vec{b}\|_{L^{n}(B_1)}\leq k(n)$, then $u(x) \geq 0$ in $B_1$.
 \end{theorem}

 \begin{remark}
 It follows from Theorem \ref{thm4} that if there exists some positive constant $k(n,p)$ such that $\|\vec{b}\|_{L^{p}(B_1)}\leq k(n,p)$ as $p\geq n$, the maximum principle for Laplacian with first order term still holds.
\end{remark}

Next, we show the strong maximum principle for Laplacian with first order term.
\begin{theorem}\label{thm5}
  Assume that $\vec{b}(x)\in L^n(B_1)$, and $u(x)\in H^1(B_1)$($n \geq 3$) is a weak solution of
   \begin{eqnarray}\label{5}
 \left\{\begin{array}{l} -\Delta u(x)+\vec{b}(x)\cdot \nabla u(x)\geq 0 \ \ \ \mbox{in} \ B_1  \\
u(x) \geq m>0 \ \ \ \mbox{on}\ \partial B_1.
\end{array}
\right.
 \end{eqnarray}
 There exists a positive constant $k(n)$ such that if
 $\|\vec{b}\|_{L^{n}(B_1)}\leq k(n)$, then $u(x) \geq m$ in $B_1$.
\end{theorem}

In the following, we consider the maximum principle and strong maximum principle for fractional Laplacian. Lots of efforts have been made on this study, see \cite{CaS, MT, QA, Si, C.L.L, C-H-L, CLL, ZL, Zli, Zli2, L.W.X, L.L.W.X, JW} and references therein.  As is well-known, the standard Laplacian in an $n$-dimensional domain possesses an explanation with regard to the diffusion, and occurs in differential equations that describe many physical phenomena. In recent years, there have been a lot of fruitful works on anomalous diffusion which is extensively observed in physics, chemistry, and biology. To characterize anomalous diffusion phenomena, the fractional Laplacian is introduced and has been widely used to model diverse
physical phenomena, such as molecular dynamics, turbulence and
water waves, and quasi-geostrophic flows \cite{CV,CM,C,TZ}. Furthermore, the fractional
Laplacian has a long history and various applications in probability and finance \cite{B1,AD}. In particular, the fractional Laplacian can be understood as the infinitesimal generator of stable L\'{e}vy process.

In comparison with Laplacian which is a local operator, the fractional Laplacian is non-local and does not act by pointwise differentiation. However, the fractional Laplacian can be defined by a global integration with respect to a singular kernel, taking the form
\begin{eqnarray*}
(-\Delta)^su(x)=C_{n,s}P.V.\int_{\mathbf{R}^{n}}\frac{u(x)-u(y)}{|x-y|^{n+2s}}dy,
\end{eqnarray*}
where $s$ is any real number between 0 and 1, and $P.V.$ stands for
the Cauchy principle value.

Note that the operator $(-\Delta)^s$ is well defined if $u(x)\in \mathcal{L}_{2s}\cap C_{loc}^{1,1}(\mathbf{R}^n)$, where
\begin{eqnarray*}
\mathcal{L}_{2s}=\left\{u: \mathbf{R}^n\rightarrow
\mathbf{R}|\int_{\mathbf{R}^n}\frac{|u(x)|}{1+|x|^{n+2s}}dx<\infty\right\}.
\end{eqnarray*}
On the other hand, if we consider $u(x)$ in the sense of distribution, only condition $u(x)\in \mathcal{L}_{2s}$ is required. More precisely, for any $\psi(x)\in C_0^{\infty}(\mathbf{R}^n)$,
\begin{eqnarray*}
(-\Delta)^su(\psi)=\int_{\mathbf{R}^n}u(x)(-\Delta)^s{\psi}(x)dx.
\end{eqnarray*}

For a long time, it is difficult to investigate the fractional Laplacian due to its
non-locality. Many scientists made a strong research effort to extend the works on Laplacian
to fractional Laplacian. Caffarelli and Silvestre \cite{CS} introduced the extension method, which transforms a
non-local problem to a local one in higher dimensions. This method has been used to show the maximum principle and strong maximum principle for fractional Laplacian by many researchers, see \cite{CaS, MT, QA} and references therein.

Here, we study the maximum principle and strong maximum principle for fractional Laplacian with critical integrability. Instead of the extension method introduced by Caffarelli and Silvestre, we work directly on the fractional Laplcian. An important issue is to study the maximum principle for fractional superharmonic functions. In \cite[Proposition 2.17]{Si},  Silvestre established the following maximum principle.

\begin{pro}\label{super}
Let $\Omega$ be an open set in $\mathbf{R}^n$. Let $u(x)\in \mathcal{L}_{2s}$ be a lower semi-continuous function in $\bar{\Omega}$ such that
\begin{eqnarray}
 \left\{\begin{array}{l} (-\Delta)^s u(x)\geq 0 \ \ \ \mbox{in} \ \Omega  \\
u(x)\geq 0 \ \ \ \mbox{in}\ {\Omega}^c
\end{array}
\right.
 \end{eqnarray}
in the sense of distribution, then $u(x)\geq 0$ in $\mathbf{R}^n$.
\end{pro}

Later, Chen, Li and Li \cite[Theorem 2.1]{C.L.L} provided a simpler proof for Proposition \ref{super} by adding the regularity condition $u(x)\in C_{loc}^{1,1}(\Omega)$. They proved the maximum principle for fractional superharmonic functions by using the integral definition of the fractional Laplcain directly.

Note that the maximum principle for fractional superharmonic functions proved in \cite{Si} and \cite{C.L.L} requires that the function is lower semi-continuous on $\bar{\Omega}$. In the next theorem, we improve the results in \cite{C.L.L, Si} by weakening the assumptions on $u$.

Throughout this paper, if a function $f(x)$ satisfies $f(x)\geq 0$ in $B_1$ in the sense of distribution, then we write
\begin{eqnarray*}
f(x)\geq 0 \ \ \text{ in } \mathcal{D}'(B_1).
\end{eqnarray*}

 \begin{theorem}\label{thm6}
 Assume that $u(x)\in \mathcal{L}_{2s}\cap L^{\frac{1}{1-s}}(B_1)$ satisfies
\begin{eqnarray}{\label{6}}
 \left\{\begin{array}{l} (-\Delta)^s u(x)\geq 0 \ \ \ \mbox{in} \ \mathcal{D}'(B_1)  \\
u(x)\geq 0 \ \ \ \mbox{in}\ B_1^c,
\end{array}
\right.
 \end{eqnarray}
then $u(x)\geq 0$ in $B_1$.
\end{theorem}
\begin{remark}
Li and Liu \cite{LLL} showed the uniqueness of $u(x)$ if $u(x)\in \mathcal{L}_{2s}\cap L^{\frac{1}{1-s}}(B_1)$ satisfies $(-\Delta)^s u(x)= 0$ in $\mathcal{D}'(B_1)$ and $u(x)= 0$ in $B_1^c$.
\end{remark}

The study on maximum principle for Schr\"{o}dinger operators involving the fractional Laplacian $(-\Delta)^s+c(x)$ has attracted a lot of attentions owing to its applications in analysis of partial differential equation, for instance, the method of moving planes for the fractional Laplacian, see \cite{C.L.L, C-H-L, CLL, ZL} and references therein. However, the condition that $c(x)$ is bounded below always requires.

In the following, we give a weaker condition for $c(x)$ to ensure the maximum principle for Schr\"{o}dinger operators involving the fractional Laplacian. More precisely, our maximum principle can be stated as follows.

 \begin{theorem}\label{thm7}
 Assume that $c(x)\in L^{\frac{n}{2s}}(B_1) $, and $u(x)\in \mathcal{L}_{2s}\cap L^{\frac{1}{1-s}}(B_1)$, $B_1\subset \mathbf{R}^n$($n\geq 3$) satisfies
   \begin{eqnarray}\label{7}
 \left\{\begin{array}{l} (-\Delta)^{s} u(x)+c(x)u(x)\geq 0 \ \ \ \mbox{in} \ \mathcal{D}'(B_1)  \\
u(x) \geq 0 \ \ \ \mbox{in}\  B_1^c.
\end{array}
\right.
 \end{eqnarray}
  There exists a positive constant $k(n,s)$ such that if
 $\|c^-\|_{L^{\frac{n}{2s}}(B_1)}\leq k(n,s)$, then $u(x) \geq 0$ in $B_1$. Here $c^{-}(x)=-\mbox{min}\{c(x),0\}$.
 \end{theorem}

The refined version of the strong maximum principle for fractional Laplacian with zero order term was discussed by Li, Wu and Xu \cite{L.W.X} which is closely related to the B$\hat{o}$cher-type theorem. They investigated the strong maximum principle on a punctured ball and required that $c(x)$ is bounded. Based on the results in Theorem \ref{thm6}, we weaken the requirements on $c(x)$.

\begin{theorem}\label{thms}
 Assume that $c(x)\in L^{p}(B_1) $$(p>\frac{n}{2s})$, and $u(x)\in \mathcal{L}_{2s}\cap L^{\frac{1}{1-s}}(B_1)$($n\geq 3$) satisfies
   \begin{eqnarray}\label{09}
 \left\{\begin{array}{l} (-\Delta)^s u(x)+c(x)u(x)\geq 0 \ \ \ \mbox{in} \ \mathcal{D}'(B_1)   \\
u(x) \geq m>0 \ \ \ \mbox{in}\  B_2\setminus B_1\\
u(x) \geq 0  \ \ \ \mbox{in} \ B_2^c.
\end{array}
\right.
 \end{eqnarray}
 There exists a positive constant $k(n,s,p)$ such that if
 $\|{c^+}\|_{L^{p}(B_1)}\leq k(n,s,p)$, then
 \begin{eqnarray*}
  u(x) \geq \gamma m \ \ \mbox{in} \ B_1,
 \end{eqnarray*}
 where $\gamma$ is a positive constant depending only on $n$, $p$ and $\|{c^+}\|_{L^p(B_1)}$.
 \end{theorem}

 Finally, with the aid of Theorem \ref{thm6}, we show the maximum principle and strong maximum principle for the fractional Laplacian with first order term.

 \begin{theorem}\label{thm8}
   Assume that $\vec{b}(x)\in L^{\frac{n}{2s}}(B_1) $, and $u(x)\in \mathcal{L}_{2s}\cap L^{\frac{1}{1-s}}(B_1)$ with $s\in (\frac{1}{2},1)$($n\geq 3$) satisfies
   \begin{eqnarray}\label{8}
 \left\{\begin{array}{l} (-\Delta)^{s} u(x)+\vec{b}(x)\cdot \nabla u(x)\geq 0 \ \ \ \mbox{in} \ \mathcal{D}'(B_1)   \\
u(x) \geq 0 \ \ \ \mbox{in}\  B_1^c.
\end{array}
\right.
 \end{eqnarray}
  There exists a positive constant $k(n,s)$ such that if
 $\|\vec{b}\|_{W^{1,\frac{n}{2s}}(B_1)}\leq k(n,s)$ and $\|\frac{\vec{b}}{d}\|_{L^{\frac{n}{2s}}(B_1)}\leq k(n,s)$, then $u(x) \geq 0$ in $B_1$. Here $d(x)=\mbox{dist}(x,\partial B_1)$.
 \end{theorem}

 \begin{theorem}\label{thm9}
   Assume that $\vec{b}(x)\in L^{\frac{n}{2s}}(B_1) $, and $u(x)\in \mathcal{L}_{2s}\cap L^{\frac{1}{1-s}}(B_1)$ with $s\in(\frac{1}{2},1)$($n\geq 3$) satisfies
   \begin{eqnarray}\label{9}
 \left\{\begin{array}{l} (-\Delta)^{s} u(x)+\vec{b}(x)\cdot \nabla u(x)\geq 0 \ \ \ \mbox{in} \ \mathcal{D}'(B_1)   \\
u(x) \geq m>0 \ \ \ \mbox{in}\  B_2\setminus B_1\\
u(x) \geq 0 \ \ \ \mbox{in}\  B_2^c.
\end{array}
\right.
 \end{eqnarray}
  There exists a positive constant $k(n,s)$ such that if
 $\|\vec{b}\|_{W^{1,\frac{n}{2s}}(B_1)}\leq k(n,s)$ and $\|\frac{\vec{b}}{d}\|_{L^{\frac{n}{2s}}(B_1)}\leq k(n,s)$, then
   \begin{eqnarray*}
  u(x) \geq \gamma m \ \ \mbox{in} \ B_1,
 \end{eqnarray*}
 where $\gamma$ is a positive constant depending only on $n$ and $s$.
 \end{theorem}

  The organization for the paper is as follows. In Section 2, we present the proofs of the different types of the maximum principle and strong maximum principle for Laplacian. The various forms of the maximum principle and strong maximum principle for fractional Laplacian are established in Section 3.

\section{Maximum principles for Laplacian}
In this section, we show the maximum principle and strong maximum principle for Laplacian with zero order term(Schr\"{o}dinger operator) and Laplacian with first order term.

\subsection{Maximum principle for Laplacian with zero order term}
In this subsection, we first show that $c(x)\in L^{\frac{n}{2}}(B_1)$ is not enough to ensure the strong maximum principle for Schr\"{o}dinger operator.
\\
\\
{\bf Proof of Theorem \ref{thm2}.} ~ Let $\alpha>0$, we define an auxiliary function
 \begin{eqnarray*}
u(x)=\frac{1}{(-ln(|x|/e))^{\alpha}},
\end{eqnarray*}
for $x\in B_1\subset \mathbf{R}^n$, {$n \geq 3$}.

Clearly, one has
  $$u(0)=0,$$
  and
  \begin{eqnarray*}
    u(x)=1\ \ \text{for } x\in \partial B_1.
   \end{eqnarray*}

A simple calculation shows that
\begin{eqnarray*}
-\Delta u(x)+c(x)u(x)=0\ \ \text{in }B_1,
\end{eqnarray*}
where
$$c(x)=\frac{\frac{\alpha(\alpha+1)}{-ln(|x|/e)}+\alpha(n-2)}{{|x|}^2(-ln(|x|/e))}.$$
Moreover, one finds that  $$c(x)\in L^{\frac{n}{2}}(B_1).$$

Now, by scaling, we define functions $u_{\epsilon}(x)$ and $c_{\epsilon}(x)$, $0<\epsilon<\frac{1}{e}$ as follows:
 \begin{eqnarray*}
u_{\epsilon}(x):=u({\epsilon}ex)=\frac{1}{(-ln({\epsilon}|x|))^{\alpha}},\ \ c_{\epsilon}(x):=c({\epsilon}ex)=\frac{\frac{\alpha(\alpha+1)}{-ln({\epsilon}|x|)}+\alpha(n-2)}{{|x|}^2(-ln({\epsilon}|x|))}.
\end{eqnarray*}
We show that $u_{\epsilon}(x)$ and $c_{\epsilon}(x)$ satisfy the required properties (i)-(iii) in Theorem \ref{thm2}.

Indeed, one has
\begin{eqnarray*}
u_{\epsilon}(0)=0,
\end{eqnarray*}
\begin{eqnarray*}
u_{\epsilon}(x)=\frac{1}{(-ln\epsilon)^{\alpha}}>0, \ \ x\in \partial B_1,
\end{eqnarray*}
and
\begin{eqnarray*}
-\Delta u_{\epsilon}(x)+c_{\epsilon}(x)u_{\epsilon}(x)=0, \ \ x\in B_1.
\end{eqnarray*}

By direct calculations, one has $c_{\epsilon}(x)\in L^{\frac{n}{2}}(B_1)$ and
$$\lim\limits_{\epsilon \rightarrow 0^+}\|c_{\epsilon}\|_{L^{\frac{n}{2}}(B_1)}=0.$$
Hence the proof of the theorem is completed.
\hfill\qedsymbol
\\

Next, we present the proof of the maximum principle for Schr\"{o}dinger operator with critical case $p=\frac{n}{2}$ in $B_1$.\\
\\
{\bf Proof of Theorem \ref{thm1}.}~ Define $u^{-}(x)=-\mbox{min}\{u(x),0\}$. It follows from (\ref{1}) that
\begin{eqnarray*}
\int_{B_1}\nabla u^{-}(x)\cdot \nabla u(x)dx+\int_{B_1}c(x)u^{-}(x)u(x)dx \geq 0.
\end{eqnarray*}
Then
\begin{eqnarray*}
\int_{B_1}|\nabla u^{-}(x)|^2dx \leq \int_{B_1}c^{-}(x)(u^{-}(x))^2dx.
\end{eqnarray*}

Clearly, $u^{-}(x)\in L^{\frac{2n}{n-2}}(B_1)$ by Sobolev embedding theorem. It follows from H\"older inequality that if
\begin{eqnarray} \label{1-0}
 \|\nabla u^{-}\|_{L^2(B_1)}\neq 0,
 \end{eqnarray}
 then
\begin{eqnarray*}
\int_{B_1}|\nabla u^{-}(x)|^2dx &\leq & \|c^-\|_{L^{\frac{n}{2}}(B_1)}\|u^{-}\|^2_{L^{\frac{2n}{n-2}}(B_1)}\\
&<& C(n)\|c^-\|_{L^{\frac{n}{2}}(B_1)}\|\nabla u^{-}\|^2_{L^2(B_1)},
\end{eqnarray*}
where $C(n)$ is a constant depending only on $n$.

Thus, one has
\begin{eqnarray*}
\|\nabla u^{-}\|_{L^2(B_1)}< (C(n)\|c^-\|_{L^{\frac{n}{2}}(B_1)})^{\frac{1}{2}}\|\nabla u^{-}\|_{L^2(B_1)}.
\end{eqnarray*}
This leads to a contradiction if
 $$\|c^-\|_{L^{\frac{n}{2}}(B_1)}\leq\frac{1}{C(n)}.$$

 Therefore, one has
 $$\|\nabla u^{-}\|_{L^2(B_1)}=0,$$
which implies $u^{-}(x)=0$ in $B_1$. Thus, one has $u(x)\geq 0$ in $B_1$.

This completes the proof of the theorem.
\hfill\qedsymbol
\\

In the following, with the aid of Theorem \ref{thm1}, we give the proof of the refined version of strong maximum principle for Schr\"{o}dinger operator if $p>\frac{n}{2}$.\\
\\
{\bf Proof of Theorem \ref{thm3}.}~
It follows from Theorem \ref{thm1} that
\begin{eqnarray*}
u(x)\geq 0 \ \ \mbox{in} \ B_1.
\end{eqnarray*}
Then one can rewrite (\ref{3}) as
\begin{eqnarray}\label{3.1}
 \left\{\begin{array}{l} -\Delta u(x)+c^{+}(x)u(x)\geq c^{-}(x)u(x)\geq 0 \ \ \ \mbox{in} \ B_1  \\
u(x) \geq m>0 \ \ \ \mbox{on}\ \partial B_1.
\end{array}
\right.
 \end{eqnarray}
Clearly, one has $c^{+}(x)\in L^p(B_1) \ (p>\frac{n}{2})$.

Suppose that $f(x)$ satisfies the following Dirichlet problem
\begin{eqnarray}\label{3.2}
 \left\{\begin{array}{l} -\Delta f(x)=c^{+}(x) \ \ \ \mbox{in} \ B_1  \\
f(x) = 0 \ \ \ \mbox{on}\ \partial B_1.
\end{array}
\right.
 \end{eqnarray}
By the classical elliptic estimates, one derives that
$$\|f\|_{L^{\infty}(B_1)}\leq C(n,p)\|c^{+}\|_{L^p(B_1)},$$
where $C(n,p)$ is a positive constant.

Let $w(x)=u(x)-m+mf(x)$. Combining (\ref{3.1}) and (\ref{3.2}) yields that $w(x)$ satisfies the following equation
\begin{eqnarray}\label{3.3}
 \left\{\begin{array}{l} -\Delta w(x)+c^{+}(x)w(x)\geq 0 \ \ \ \mbox{in} \ B_1  \\
w(x) \geq 0 \ \ \ \mbox{on}\ \partial B_1.
\end{array}
\right.
 \end{eqnarray}
It follows from Theorem \ref{thm1} that $w(x)\geq 0$ in $B_1$.

Therefore, one has, for $x\in B_1$,
\begin{eqnarray*}
u(x)&\geq& m(1-f(x))\\
&\geq& m(1-\|f\|_{L^{\infty}(B_1)})\\
&\geq &m(1-C(n,p)\|c^{+}\|_{L^p(B_1)}).
\end{eqnarray*}
Let
$$k(n,p)=\frac{1}{2C(n,p)}.$$

If $\|c^{+}\|_{L^p(B_1)}\leq k(n,p)$, then one has
$$u(x)\geq \frac{m}{2}.$$

Hence the proof of the theorem is completed.
\hfill\qedsymbol
\\

\subsection{Maximum principles for Laplacian with first order term}
In this subsection, we first give the proof of the maximum principle for Laplacian with first order term.\\
\\
{\bf Proof of Theorem \ref{thm4}.}
\quad ~ Define $u^{-}(x)=-\mbox{min}\{u(x),0\}$. It follows from (\ref{4}) that
\begin{eqnarray*}
\int_{B_1}\nabla u^{-}(x)\cdot \nabla u(x)dx+\int_{B_1}(\vec{b}(x)\cdot \nabla u(x))u^{-}(x)dx \geq 0.
\end{eqnarray*}
Then
\begin{eqnarray*}
\int_{B_1}|\nabla u^{-}(x)|^2dx \leq \int_{B_1}|\vec{b}(x)||\nabla u^{-}(x)||u^{-}(x)|dx.
\end{eqnarray*}
Since $u^{-}(x)\in H^1(B_1)$, using Sobolev embedding theorem, one can obtain $u^{-}(x)\in L^{\frac{2n}{n-2}}(B_1)$.

 It follows from H\"older inequality that if $\|\nabla u^{-}\|_{L^2(B_1)}\neq 0$, one has
\begin{eqnarray*}
\int_{B_1}|\nabla u^{-}(x)|^2dx &\leq & \|\vec{b}\|_{{L^n}(B_1)}\|\nabla u^{-}\|_{{L^2}(B_1)}\|u^{-}\|_{L^{\frac{2n}{n-2}}(B_1)}\\
&<& C(n)\|\vec{b}\|_{{L^n}(B_1)}\|\nabla u^{-}\|^2_{L^2(B_1)}.
\end{eqnarray*}
Then
\begin{eqnarray*}
\|\nabla u^{-}\|_{L^2(B_1)}< (C(n)\|\vec{b}\|_{{L^n}(B_1)})^{\frac{1}{2}}\|\nabla u^{-}\|_{L^2(B_1)},
\end{eqnarray*}
where $C(n)$ is a constant depending only on $n$.

This leads to a contradiction if
 $$\|\vec{b}\|_{{L^n}(B_1)}\leq\frac{1}{C(n)}.$$
Therefore, one has
 $$\|\nabla u^{-}\|_{L^2(B_1)}=0.$$
This implies that $u^{-}(x)=0$ in $B_1$. Thus, one has $u(x)\geq 0$ in $B_1$.

Hence the proof of the theorem is completed.
\hfill\qedsymbol
\\

Next, we show the strong maximum principle for Laplacian with first order term.\\
\\
{\bf Proof of Theorem \ref{thm5}}
\quad ~ Let $v(x)=u(x)-m$. One can derive that
\begin{eqnarray}\label{5.1}
 \left\{\begin{array}{l} -\Delta v(x)+\vec{b}(x)\cdot \nabla v(x)\geq 0 \ \ \ \mbox{in} \ B_1  \\
v(x) \geq 0 \ \ \ \mbox{on}\ \partial B_1.
\end{array}
\right.
 \end{eqnarray}
It follows from Theorem \ref{thm4} that $v(x) \geq 0$ in $B_1$. Thus,
\begin{eqnarray*}
u(x)\geq m \ \ \mbox{in}\ B_1.
\end{eqnarray*}

This completes the proof of the theorem.

\section{Maximum principles for fractional Laplacian}
\quad ~ The aim of this section is to prove the maximum principle and strong maximum principle for fractional Laplacian. These results are useful in analysis of the fractional Laplace equations and play an important role in understanding the nonlocal property for the fractional Laplacian.
\subsection{Preliminaries}
\quad ~In this subsection, we first introduce the explicit formula for the Poisson kernel and Green function for the fractional Laplacian. Then we present two basic lemmas which are helpful to prove the maximum principles for the fractional Laplacian.

The formula for the Poisson kernel of balls was obtained by Riesz in \cite{Rie}, and the formula for Green function of balls was obtained by Blumenthal, Getoor and Ray in \cite{BGR}. Specifically, let $r>0$. For any $x\in B_r$ and $y\in \bar{B}_r^c$, the Poisson kernel $P_r(x,y)$ is defined by
\begin{eqnarray*}
 P_r(x,y)=c(n,s)\left(\frac{r^2-|x|^2}{|y|^2-r^2}\right)^s\frac{1}{|x-y|^n},
 \end{eqnarray*}
where $c(n,s)$ is a constant depending only on $n$ and $s$.
For any $x,z\in B_r$ and $x\neq z$, the Green function $G_r(x,z)$ is given by
\begin{eqnarray*}
 G_r(x,z)=\Phi(x-z)-\int_{B_r^c}\Phi(z-y)P_r(x,y)dy,
 \end{eqnarray*}
 where
 \begin{eqnarray*}
 \Phi(x-z)=\frac{c(n,s)}{|x-z|^{n-2s}}.
 \end{eqnarray*}

 Moreover, the explicit formula for the Green function is that for fixed $r>0$, $n> 2s$,
\begin{eqnarray}\label{green}
G_r(x,z)=\kappa(n,s)|z-x|^{2s-n}\int_{0}^{\sigma_r(x,z)}\frac{t^{s-1}}{(t+1)^{\frac{n}{2}}}dt,
\end{eqnarray}
where
\begin{eqnarray*}
\sigma_r(x,z)=\frac{(r^2-|x|^2)(r^2-|z|^2)}{r^2|x-z|^2}
\end{eqnarray*}
and $\kappa(n,s)$ is a constant depending only on $n$ and $s$.

Bucur \cite{C.B} showed that if $g(x)\in \mathcal{L}_{2s}\cap C(\mathbf{R}^n)$, then
\begin{eqnarray}\label{6.1.1}
 u_g(x)=\left\{\begin{array}{l} \int_{B_r^c}P_r(x,y)g(y)dy \ \ \ \mbox{if} \ x\in B_r  \\
g(x)  \ \ \ \mbox{if}\ x\in B_r^c
\end{array}
\right.
 \end{eqnarray}
 is the unique pointwise continuous solution of the equation
 \begin{eqnarray}\label{6.1.2}
 \left\{\begin{array}{l} (-\Delta)^s u(x)=0 \ \ \ \mbox{in} \ B_r  \\
u(x) = g(x)  \ \ \ \mbox{in}\ B_r^c.
\end{array}
\right.
 \end{eqnarray}

On the other hand, if $h(x)\in C^{2s+\epsilon}(B_r)\cap C(\bar{B_r})$, then
\begin{eqnarray}\label{6.1.3}
 u(x)=\left\{\begin{array}{l} \int_{B_r}h(y)G_r(x,y)dy \ \ \ \mbox{if} \ x\in B_r  \\
0  \ \ \ \mbox{if}\ x\in B_r^c
\end{array}
\right.
 \end{eqnarray}
 is the unique pointwise continuous solution of the equation
 \begin{eqnarray}\label{6.1.4}
 \left\{\begin{array}{l} (-\Delta)^s u(x)=h(x) \ \ \ \mbox{in} \ B_r  \\
u(x) = 0  \ \ \ \mbox{in}\ B_r^c.
\end{array}
\right.
 \end{eqnarray}


Now, we give two basic lemmas for the fractional Laplacian.

\begin{lemma}{\label{lemma1.2}}
Let $\eta_{\epsilon}(x)=\frac{1}{{\epsilon}^n}\eta(\frac{x}{\epsilon})$, where $\eta(x)$ satisfies  $\eta(x) \in C_0^{\infty}(B_1)$, $\eta(x)\geq 0$ and $\int_{B_1}\eta(x) dx=1$. Assume that $u(x) \in \mathcal{L}_{2s}$ satisfies
\begin{eqnarray}
(-\Delta)^s u(x)\leq 0 \ \ \mbox{in}\ \mathcal{D}'(B_1),
\end{eqnarray}
and define its mollification $u_{\epsilon}(x)= \eta_{\epsilon}\ast u(x)$ in $B_{1-\epsilon}$. Then $u_{\epsilon}(x)$ satisfies
\begin{eqnarray}
(-\Delta)^s u_{\epsilon}(x)\leq 0 \ \ \mbox{in}\ B_{1-\epsilon}.
\end{eqnarray}
\end{lemma}
\begin{proof}
Taking into account the definition of the fractional Laplacian and mollification, it follows that for $x\in B_{1-\epsilon}$, one has
\begin{eqnarray*}
(-\Delta)^su_{\epsilon}(x)&=&C_{n,s}P.V.\int_{\mathbf{R}^{n}}\frac{u_{\epsilon}(x)-u_{\epsilon}(y)}{|x-y|^{n+2s}}dy\\
                          &=&C_{n,s}P.V.\int_{\mathbf{R}^{n}}\frac{\int_{\mathbf{R}^{n}}\eta_{\epsilon}(x-z)u(z)dz-\int_{\mathbf{R}^{n}}\eta_{\epsilon}(y-z)u(z)dz}{|x-y|^{n+2s}}dy\\
                          &=&\int_{\mathbf{R}^{n}}u(z)(-\Delta)^s_x\eta_{\epsilon}(x-z)dz\\
                          &=&\int_{\mathbf{R}^{n}}u(z)(-\Delta)^s_z\eta_{\epsilon}(x-z)dz\\
                          &\leq&0.
\end{eqnarray*}
This finishes the proof of the lemma.
\end{proof}

Inspired by the idea of \cite{L.L.W.X} and \cite{L.W.X}, we give the following lemma. However, the statement in \cite[Theorem 5.4]{L.L.W.X} requires
\begin{eqnarray}
\|\vec{b}(x)\|_{C^1{(\Omega})}+\|c(x)\|_{L^{\infty}(\Omega)}<\infty,
\end{eqnarray}
where $\Omega$ is a domain in $\mathbf{R}^n$. Here, we weaken the conditions for $\vec{b}(x)$ and $c(x)$.
\begin{lemma}\label{lemma1}
 Assume that $u(x) \in \mathcal{L}_{2s}\cap L^{\frac{1}{1-s}}(B_1)$, $\vec{b}\in W^{1,\frac{1}{s}}(B_1)$ and $c(x)\in L^{\frac{1}{s}}(B_1)$ with $s\in(\frac{1}{2},1)$ satisfy
 \begin{equation}
 (-\Delta)^s u(x)+\vec{b}(x)\cdot \nabla u(x)+c(x)u(x)\geq 0\ \ \mbox{in}\  \mathcal{D}'(B_1).
\end{equation}
Then for $v(x)=min\{u(x),0\}$, one has
 \begin{equation}\label{com}
 (-\Delta)^s v(x)+\vec{b}(x)\cdot \nabla v(x)+c(x)v(x)\geq 0\ \ \mbox{in}\  \mathcal{D}'(B_1) .
\end{equation}
\end{lemma}

\begin{proof}
The proof is divided into three steps.

\textbf{Step 1. } Set $U=\{x\in B_1| u(x)<0\}$. In this step, we assume that $u(x)$ is a smooth function and $\partial U$ is $C^1$. Let $\varphi(x)\in C^{\infty}_0(B_1)$ be a nonnegative function and $\bf{n}$ the outward unit normal vector of $\partial U$. From the definition of the fractional Laplacian, one can derive that
\begin{eqnarray*}
&&\int_{\mathbf{R}^n}v(x)(-\Delta)^s\varphi(x)dx\\
&&=\int_{\mathbf{R}^n}v(x)\left(C_{n,s}\lim\limits_{\epsilon\rightarrow 0}\int_{\mathbf{R}^n\backslash B_{\epsilon}(x)}\frac{\varphi(x)-\varphi(y)}{|x-y|^{n+2s}}dy\right)dx\\
&&=\int_{U}u(x)\left(C_{n,s}\lim\limits_{\epsilon\rightarrow 0}\int_{\mathbf{R}^n\backslash B_{\epsilon}(x)}\frac{\varphi(x)-\varphi(y)}{|x-y|^{n+2s}}dy\right)dx\\
&&=\int_{U}C_{n,s}\lim\limits_{\epsilon\rightarrow 0}\int_{\mathbf{R}^n\backslash B_{\epsilon}(x)}\frac{(u(x)-u(y))\varphi(x)}{|x-y|^{n+2s}}dydx\\
&&~~~~~~~~~~~~~~~~~~~~~~~~+\int_{U}C_{n,s}\lim\limits_{\epsilon\rightarrow 0}\int_{\mathbf{R}^n \backslash B_{\epsilon}(x)}\frac{u(y)\varphi(x)-u(x)\varphi(y)}{|x-y|^{n+2s}}dydx\\
&&=\int_{U}(-\Delta)^su(x)\varphi(x)dx+\int_{U}C_{n,s}\lim\limits_{\epsilon\rightarrow 0}\int_{U\backslash B_{\epsilon}(x)}\frac{u(y)\varphi(x)-u(x)\varphi(y)}{|x-y|^{n+2s}}dydx\\
&&~~~~~~~~~~~~~~~~~~~~~~~~~~~~~~~~~~~~~+\int_{U}C_{n,s}\lim\limits_{\epsilon\rightarrow 0}\int_{U^c\backslash B_{\epsilon}(x)}\frac{u(y)\varphi(x)-u(x)\varphi(y)}{|x-y|^{n+2s}}dydx.
\end{eqnarray*}

Clearly, one has
\begin{eqnarray*}
\int_{U}C_{n,s}\lim\limits_{\epsilon\rightarrow 0}\int_{U\backslash B_{\epsilon}(x)}\frac{u(y)\varphi(x)-u(x)\varphi(y)}{|x-y|^{n+2s}}dydx=0,
\end{eqnarray*}
and
\begin{eqnarray*}
\int_{U}C_{n,s}\lim\limits_{\epsilon\rightarrow 0}\int_{U^c\backslash B_{\epsilon}(x)}\frac{u(y)\varphi(x)-u(x)\varphi(y)}{|x-y|^{n+2s}}dydx\geq 0.
\end{eqnarray*}
Thus, one has
\begin{eqnarray}\label{1-1}
\int_{\mathbf{R}^n}v(x)(-\Delta)^s\varphi(x)dx\geq \int_{U}(-\Delta)^su(x)\varphi(x)dx.
\end{eqnarray}

On the other hand, for the lower order term, one has
\begin{eqnarray}\label{1-2}
&&\int_{B_1}v(x)\left(-div(\vec{b}(x)\varphi(x))+c(x)\varphi(x)\right)\\ \nonumber
&&=\int_{U}u(x)\left(-div(\vec{b}(x)\varphi(x))+c(x)\varphi(x)\right)\\  \nonumber
&&=\int_{U}(\vec{b}(x)\cdot \nabla u(x)+c(x)u(x))\varphi(x)dx-\int_{\partial U}u(x)\varphi(x)\vec{b}(x)\cdot \textbf{n}\ dS\\
&&=\int_{U}(\vec{b}(x)\cdot \nabla u(x)+c(x)u(x))\varphi(x)dx. \nonumber
\end{eqnarray}
Combining (\ref{1-1}) and (\ref{1-2}) leads to
\begin{eqnarray*}
(-\Delta)^s v(x)+\vec{b}(x)\cdot \nabla v(x)+c(x)v(x)\geq \left((-\Delta)^s u(x)+\vec{b}(x)\cdot \nabla u(x)+c(x)u(x)\right)\chi_{\{u(x)<0\}} \ \ \mbox{in} \ \mathcal{D}'(B_1).
\end{eqnarray*}

Thus, one can obtain that
\begin{eqnarray}\label{inequality}
(-\Delta)^s v(x)+\vec{b}(x)\cdot \nabla v(x)+c(x)v(x)\geq 0  \ \ \mbox{in} \ \mathcal{D}'(B_1) .
\end{eqnarray}

\textbf{Step 2. } In this step, we assume that $u(x)$ is a smooth function, however $\partial U$ may not be $C^1$. We show that (\ref{inequality}) still holds.

Let $u_\theta(x)=u(x)-\theta$ for $\theta<0$. By Sard's theorem, one can choose a non-positive monotone increasing sequence $\{\theta_j\}^{\infty}_{j=1}$ satisfying $\theta_j\rightarrow 0$, such that the set $\tilde{U}=\{x\in B_1|u_{\theta_j}(x)<0\}$ satisfies $\partial \tilde{U}\in C^1$ for each $j$.

Denote $v_{\theta_j}(x)=min\{u_{\theta_j}(x),0\}$. It follows from the results in Step 1 that
\begin{eqnarray}
&&(-\Delta)^s v_{\theta_j}(x)+\vec{b}\cdot \nabla v_{\theta_j}(x)+c(x) v_{\theta_j}(x)\\
&&\geq \left((-\Delta)^s u_{\theta_j}(x)+\vec{b}\cdot \nabla u_{\theta_j}(x)+c(x) u_{\theta_j}(x)\right)\chi_{\{u(x)<\theta_j\}} \geq -\theta_jc(x)\chi_{\{u(x)<\theta_j\}} \ \ \mbox{in} \ \mathcal{D}'(B_1).\nonumber
\end{eqnarray}

Thus, for $\varphi(x)\in C^{\infty}_0(B_1)$, one derives that
\begin{eqnarray}\label{1-3}
\int_{B_1}v_{\theta_j}(x)\left((-\Delta)^s\varphi(x)-div(\vec{b}(x)\varphi(x))+c(x)\varphi(x)\right)dx\geq \int_{B_1} -\theta_jc(x)\varphi(x)\chi_{\{u(x)<\theta_j\}}dx.
\end{eqnarray}

Clearly, one has
 $$\bigcup_{j=1}^{\infty}\{x\in B_1|u(x)<\theta_j\}=\{x\in B_1|u(x)<0\},$$
 where $\{x\in B_1|u(x)<\theta_j\}$ is a sequence of monotonically increasing sets. One can obtain that
 \begin{eqnarray*}
 \int_{B_1}c(x)\varphi(x)\chi_{\{u(x)<\theta_j\}}dx &\leq&  \left(\int_{\{u(x)<\theta_j\}}|c(x)|^{\frac{1}{s}}dx\right)^{s}\left(\int_{\{u(x)<\theta_j\}}|\varphi(x)|^{\frac{1}{1-s}}dx\right)^{1-s}\\
 &\leq& C(n,s).
  \end{eqnarray*}
 Thus, letting $j\rightarrow \infty$ in (\ref{1-3}), one has
 \begin{eqnarray*}
 \int_{B_1}v(x)\left((-\Delta)^s\varphi(x)-div(\vec{b}(x)\varphi(x))+c(x)\varphi(x)\right)dx\geq 0 \ \ \mbox{in} \ B_1.
 \end{eqnarray*}
 This implies
 \begin{eqnarray*}
(-\Delta)^s v(x)+\vec{b}(x)\cdot \nabla v(x)+c(x)v(x)\geq 0  \ \ \mbox{in} \ \mathcal{D}'(B_1) .
\end{eqnarray*}

\textbf{Step 3.} In the last step, we assert that (\ref{inequality}) still holds even if $u(x)$ may not be a smooth function and $\partial U$ may not be $C^1$. Actually, one can derive that its mollification $u_{\epsilon}(x)$ satisfies
\begin{eqnarray*}
(-\Delta )^s u_{\epsilon}(x)+\vec{b}(x)\cdot \nabla u_{\epsilon}(x)+c(x)u_{\epsilon}(x) \geq P_{\epsilon}(x)+Q_{\epsilon}(x)\ \ \mbox{in}\ B_{1-\epsilon},
\end{eqnarray*}
where $P_{\epsilon}(x)=\vec{b}(x)\cdot \nabla u_{\epsilon}(x)- (\vec{b}\cdot \nabla u)_{\epsilon}(x)$ and $Q_{\epsilon}(x)=c(x)u_{\epsilon}(x)-(cu)_{\epsilon}(x)$.

It follows from the results in Step 2 that
 \begin{eqnarray}\label{l-1}
(-\Delta)^s v_{\epsilon}(x)+\vec{b}(x)\cdot \nabla v_{\epsilon}(x)+c(x)v_{\epsilon}(x)\geq (P_{\epsilon}(x)+Q_{\epsilon}(x))\chi_{\{u_{\epsilon}(x)<0\}}\ \ \mbox{in}\ \mathcal{D}'(B_{1-\epsilon}),
\end{eqnarray}
where $v_{\epsilon}(x)=min\{u_{\epsilon}(x),0\}$. Once we know that $P_{\epsilon}(x)\rightarrow 0$ and $Q_{\epsilon}(x)\rightarrow 0$ in $L^1_{loc}(B_1)$ as $\epsilon \rightarrow 0^{+}$, then by letting $\epsilon \rightarrow 0^{+}$ in (\ref{l-1}), one can obtain
\begin{eqnarray*}
(-\Delta)^s v(x)+\vec{b}(x)\cdot \nabla v(x)+c(x)v(x)\geq 0  \ \ \mbox{in} \ \mathcal{D}'(B_1)
\end{eqnarray*}
as desired.

In the following, we prove $P_{\epsilon}(x)\rightarrow 0$ and $Q_{\epsilon}(x)\rightarrow 0$ in $L^1_{loc}(B_1)$ as $\epsilon \rightarrow 0^{+}$, respectively.

(i) Prove $P_{\epsilon}(x)\rightarrow 0$ in $L^1_{loc}(B_1)$ as $\epsilon \rightarrow 0^{+}$.
\begin{eqnarray*}
P_{\epsilon}(x)&=&\vec{b}(x)\cdot \nabla u_{\epsilon}(x)- (\vec{b}\cdot \nabla u)_{\epsilon}(x)\\
               &=&\int_{B_{\epsilon}(x)}u(y)\vec{b}(x)\cdot \nabla_{x}\eta_{\epsilon}(x-y)dy+\int_{B_{\epsilon}(x)}u(y)div_y(\vec{b}(y)\eta_{\epsilon}(x-y))dy\\
               &=&\int_{B_{\epsilon}(x)}u(y)div_y[(\vec{b}(y)-\vec{b}(x))\eta_{\epsilon}(x-y)]dy\\
               &=&\int_{B_{\epsilon}(x)}(u(y)-u(x))div_y[(\vec{b}(y)-\vec{b}(x))\eta_{\epsilon}(x-y)]dy\\
               &=&\int_{B_{\epsilon}(x)}(u(y)-u(x))div_y\vec{b}(y)\eta_{\epsilon}(x-y)dy+\int_{B_{\epsilon}(x)}(u(y)-u(x))(\vec{b}(y)-\vec{b}(x))\cdot \nabla_y\eta_{\epsilon}(x-y)dy\\
               &:=&I_1(x)+I_2(x).
\end{eqnarray*}
Then, one has
\begin{eqnarray*}
\int_{B_{1-\epsilon}}|I_1(x)|dx &\leq& \int_{B_{1-\epsilon}}\int_{B_{\epsilon}(x)}|u(y)-u(x)||div_y\vec{b}(y)|\eta_{\epsilon}(x-y)dydx\\
                                &=& \int_{B_{1-\epsilon}}\int_{B_{\epsilon}(0)}|u(x+z)-u(x)||div_z\vec{b}(x+z)|\eta_{\epsilon}(z)dzdx\\
                                &=& \int_{B_{\epsilon}(0)}\left(\int_{B_{1-\epsilon}}|u(x+z)-u(x)||div_z\vec{b}(x+z)|dx\right) \eta_{\epsilon}(z) dz\\
                                &\leq&\int_{B_{\epsilon}(0)}\left((\int_{B_{1-\epsilon}}|u(x+z)-u(x)|^{\frac{1}{1-s}}dx)^{1-s}(\int_{B_{1-\epsilon}}|div_z\vec{b}(x+z)|^{\frac{1}{s}}dx)^s\right)\eta_{\epsilon}(z)dz,
\end{eqnarray*}
where we have used H\"older inequality.

It follows from the condition $\vec{b}(x) \in W^{1,\frac{1}{s}}(B_1)$ and the properties of the mollifier that
\begin{eqnarray}\label{I1}
&&\lim\limits_{\epsilon\rightarrow 0^{+}}\int_{B_{1-\epsilon}}|I_1(x)|dx \\ \nonumber
&&\leq C(n,s)\lim\limits_{\epsilon\rightarrow 0^{+}}\int_{B_{\epsilon}(0)}\left((\int_{B_{1-\epsilon}}|u(x+z)-u(x)|^{\frac{1}{1-s}}dx)^{1-s}\right)\eta_{\epsilon}(z)dz=0.
\end{eqnarray}

On the other hand, one can deduce that
\begin{eqnarray*}
\int_{B_{1-\epsilon}}|I_2(x)|dx &\leq& \int_{B_{1-\epsilon}}\int_{B_{\epsilon}(x)}|u(y)-u(x)||\vec{b}(y)-\vec{b}(x)||\nabla_y\eta_{\epsilon}(x-y)|dydx\\
                             &=& \int_{B_{1-\epsilon}}\int_{B_{\epsilon}(0)}|u(x+z)-u(x)||\vec{b}(x+z)-\vec{b}(x)||\nabla\eta_{\epsilon}(z)|dzdx\\
                             &=&\int_{B_{\epsilon}(0)}\left(\int_{B_{1-\epsilon}}|u(x+z)-u(x)||\vec{b}(x+z)-\vec{b}(x)|dx\right)|\nabla\eta_{\epsilon}(z)|dz\\
                             &=& \int_{B_{\epsilon}(0)}\left(\int_{B_{1-\epsilon}}|u(x+z)-u(x)||\int_0^1\frac{d}{dt}\vec{b}(x+tz)dt|dx\right)|\nabla\eta_{\epsilon}(z)|dz\\
                             &\leq& \int_0^1\int_{B_{\epsilon}(0)}\left(\int_{B_{1-\epsilon}}|u(x+z)-u(x)||\nabla\cdot\vec{b}(x+tz)||z|dx\right)|\nabla\eta_{\epsilon}(z)|dzdt.
\end{eqnarray*}
It follows from H\"older inequality and $\vec{b}(x) \in W^{1,\frac{1}{s}}(B_1)$ that
\begin{eqnarray*}
&&\int_{B_{1-\epsilon}}|I_2(x)|dx\\
&&\leq \int_0^1\epsilon\int_{B_{\epsilon}(0)}\left((\int_{B_{1-\epsilon}}|u(x+z)-u(x)|^{\frac{1}{1-s}}dx)^{1-s}(\int_{B_{1-\epsilon}}|\nabla\cdot\vec{b}(x+tz)|^{\frac{1}{s}}dx)^s\right)|\nabla\eta_{\epsilon}(z)|dzdt\\
&&\leq \epsilon C(n,s)\int_{B_{\epsilon}(0)}\left((\int_{B_{1-\epsilon}}|u(x+z)-u(x)|^{\frac{1}{1-s}}dx)^{1-s}\right)|\nabla\eta_{\epsilon}(z)|dz.
\end{eqnarray*}
Clearly, for $z\in B_{\epsilon}(0)$,
\begin{eqnarray*}
\epsilon\int_{B_{\epsilon}(0)}|\nabla\eta_{\epsilon}(z)|dz\leq \epsilon\int_{B_{\epsilon}(0)}\frac{C}{{\epsilon}^{n+1}}dz\leq C(n)
\end{eqnarray*}

Thus, one has
\begin{eqnarray}\label{I2}
\lim\limits_{\epsilon\rightarrow 0^{+}}\int_{B_{1-\epsilon}}|I_2(x)|dx=0.
\end{eqnarray}
Combining (\ref{I1}) and (\ref{I2}) leads to
\begin{eqnarray*}
\lim\limits_{\epsilon\rightarrow 0^{+}}\int_{B_{1-\epsilon}}|P_{\epsilon}(x)|dx=0.
\end{eqnarray*}
This implies $P_{\epsilon}(x)\rightarrow 0$ in $L^1_{loc}(B_1)$ as $\epsilon \rightarrow 0^{+}$.

(ii) Prove $Q_{\epsilon}(x)\rightarrow 0$ in $L^1_{loc}(B_1)$ as $\epsilon \rightarrow 0^{+}$.
\begin{eqnarray*}
Q_{\epsilon}(x)&=&c(x)u_{\epsilon}(x)-(cu)_{\epsilon}(x)\\
               &=&\int_{B_{\epsilon}(x)}c(x)u(y)\eta_{\epsilon}(x-y)dy-\int_{B_{\epsilon}(x)}c(y)u(y)\eta_{\epsilon}(x-y)dy\\
               &=&\int_{B_{\epsilon}(x)}u(y)(c(x)-c(y))\eta_{\epsilon}(x-y)dy.
\end{eqnarray*}
It follows that
\begin{eqnarray*}
\int_{B_{1-\epsilon}}|Q_{\epsilon}(x)|dx &\leq& \int_{B_{1-\epsilon}}\int_{B_{\epsilon}(x)}|u(y)||c(x)-c(y)|\eta_{\epsilon}(x-y)dydx\\
                                &=& \int_{B_{1-\epsilon}}\int_{B_{\epsilon}(0)}|u(x+z)||c(x+z)-c(x)|\eta_{\epsilon}(z)dzdx\\
                                &=& \int_{B_{\epsilon}(0)}\left(\int_{B_{1-\epsilon}}|u(x+z)||c(x+z)-c(x)|dx\right) \eta_{\epsilon}(z) dz\\
                                &\leq&\int_{B_{\epsilon}(0)}\left((\int_{B_{1-\epsilon}}|u(x+z)|^{\frac{1}{1-s}}dx)^{1-s}(\int_{B_{1-\epsilon}}|c(x+z)-c(x)|^{\frac{1}{s}}dx)^s\right)\eta_{\epsilon}(z)dz.
\end{eqnarray*}
It follows from the condition $c(x) \in L^{\frac{1}{s}}(B_1)$ and the properties of the mollifier that
\begin{eqnarray*}
\lim\limits_{\epsilon\rightarrow 0^{+}}\int_{B_{1-\epsilon}}|Q_{\epsilon}(x)|dx=0.
\end{eqnarray*}
This implies $Q_{\epsilon}(x)\rightarrow 0$ in $L^1_{loc}(B_1)$ as $\epsilon \rightarrow 0^{+}$.

Thus, one has
\begin{eqnarray*}
(-\Delta)^s v(x)+\vec{b}(x)\cdot \nabla v(x)+c(x)v(x)\geq 0  \ \ \mbox{in} \ \mathcal{D}'(B_1) .
\end{eqnarray*}
This completes the proof of the lemma.
\end{proof}

\subsection{Maximum principle for fractional superharmonic functions}
In this subsection, we prove the maximum principle for fractional superharmonic functions. Our strategy is to use the properties for the mollification of the fractional superharmonic function and the representation formula for the fractional Laplacian.

For convenience, we state Theorem \ref{thm6} again here.

\begin{theorem}
 Assume that $u(x)\in \mathcal{L}_{2s}\cap L^{\frac{1}{1-s}}(B_1)$ satisfies
\begin{eqnarray}{\label{3.3}}
 \left\{\begin{array}{l} (-\Delta)^s u(x)\geq 0 \ \ \ \mbox{in} \ \mathcal{D}'(B_1)  \\
u(x)\geq 0 \ \ \ \mbox{in}\ B_1^c,
\end{array}
\right.
 \end{eqnarray}
then $u(x)\geq 0$ in $B_1$.
\end{theorem}
\begin{proof}
Set $u^-(x)=-\text{min}\{u(x),0\}$. It follows from Lemma \ref{lemma1} that
\begin{eqnarray}{\label{6-1}}
 \left\{\begin{array}{l} (-\Delta)^s u^{-}(x)\leq 0 \ \ \ \mbox{in} \ \mathcal{D}'(B_1)  \\
u^{-}(x)= 0 \ \ \ \mbox{in}\ B_1^c.
\end{array}
\right.
 \end{eqnarray}

Then for $\epsilon>0$, it follows from Lemma \ref{lemma1.2} that the mollification $u_{\epsilon}^{-}(x)$ satisfies
\begin{eqnarray}{\label{6-2}}
 \left\{\begin{array}{l} (-\Delta)^s u_{\epsilon}^{-}(x)=f(x)\leq 0 \ \ \ \mbox{in} \ B_{1-\epsilon}  \\
u_{\epsilon}^{-}(x)= 0 \ \ \ \mbox{in}\ B_{1+\epsilon}^c,
\end{array}
\right.
 \end{eqnarray}
 where $f(x)$ is a smooth function.

 For any $r\in (0, 1-\epsilon$], applying the representation formula of $u_{\epsilon}^{-}(x)$, one has
 \begin{eqnarray}\label{6-3}
 u_{\epsilon}^{-}(x)&=&\int_{B_r^c}P_r(x,y)u_{\epsilon}^{-}(y)dy+\int_{B_r}f(y)G_r(x,y)dy\\\nonumber
                    &\leq&\int_{B_r^c}P_r(x,y)u_{\epsilon}^{-}(y)dy.
 \end{eqnarray}

  Now, for fixed $t\in (0,1)$ and any $x\in B_t$, one can choose $\epsilon\in (0, \frac{1-t}{3})$, take average of the right side of (\ref{6-3}) with $r\in [1-2\epsilon, 1-\epsilon]$ and deduce that
 \begin{eqnarray*}
 u_{\epsilon}^{-}(x)&\leq&\frac{1}{\epsilon}\int_{1-2\epsilon}^{1-\epsilon}dr\int_{B_r^c}P_r(x,y)u_{\epsilon}^{-}(y)dy\\
                    &=&\frac{1}{\epsilon}\int_{B_{1-2\epsilon}^c}\int_{1-2\epsilon}^{min\{1-\epsilon,|y|\}}P_r(x,y)u_{\epsilon}^{-}(y)drdy\\
                    &=&\frac{1}{\epsilon}\int_{B_{1-2\epsilon}^c}\int_{1-2\epsilon}^{min\{1-\epsilon,|y|\}}\left(\frac{r^2-|x|^2}{|y|^2-r^2}\right)^s\frac{c(n,s)u_{\epsilon}^{-}(y)}{|x-y|^n}drdy.
 \end{eqnarray*}
 Clearly, for $x\in B_t$ and $y\in B^c_{1-2\epsilon}$,
 \begin{eqnarray}\label{0-3}
 |x-y|^n\geq (|y|-|x|)^n\geq(1-2\epsilon-t)^n\geq \left(\frac{1-t}{3}\right)^n,
 \end{eqnarray}
 and
 \begin{eqnarray}\label{0-2}
 \frac{1}{(|y|^2-r^2)^s}=\frac{1}{(|y|-r)^s}\frac{1}{(|y|+r)^s}\leq \frac{C(t)}{(|y|-r)^s}.
 \end{eqnarray}

 Combining (\ref{0-3}) and (\ref{0-2}) yields
  \begin{eqnarray*}
 &&\int_{1-2\epsilon}^{min\{1-\epsilon,|y|\}}\left(\frac{r^2-|x|^2}{|y|^2-r^2}\right)^s\frac{c(n,s)}{|x-y|^n}dr\\
 &&\leq C(t,n,s)\int_{1-2\epsilon}^{min\{1-\epsilon,|y|\}}\frac{1}{(|y|-r)^s}dr\\
 &&\leq \frac{ C_1(t,n,s)}{\epsilon^{s-1}}.
 \end{eqnarray*}
It follows that
 \begin{eqnarray*}
u_{\epsilon}^{-}(x) &\leq&\frac{C_1(t,n,s)}{{\epsilon}^s}\int_{ B_{1-2\epsilon}^c}u_{\epsilon}^{-}(y)dy\\
&\leq&\frac{C_2(t,n,s)}{{\epsilon}^s}\int_{B_1\backslash B_{1-3\epsilon}}u^{-}(y)dy\\
                    &\leq&\frac{C_2(t,n,s)}{{\epsilon}^s}|B_1\backslash B_{1-3\epsilon}|^s\left(\int_{B_1\backslash B_{1-3\epsilon}}(u^{-}(y))^{\frac{1}{1-s}}dy\right)^{1-s}\\
                    &\leq&C_3(t,n,s)\left(\int_{B_1\backslash B_{1-3\epsilon}}(u^{-}(y))^{\frac{1}{1-s}}dy\right)^{1-s}.
 \end{eqnarray*}

Clearly,  $u^{-}(x)\in L^{\frac{1}{1-s}}(B_1)$ and
\begin{eqnarray*}
C_3(t,n,s)\left(\int_{B_1\backslash B_{1-3\epsilon}}(u^{-}(y))^{\frac{1}{1-s}}dy\right)^{1-s}\rightarrow 0 \ \ \mbox{as} \ \epsilon\rightarrow 0^{+}.
\end{eqnarray*}
Hence, one has
\begin{eqnarray*}
\|u^-\|_{{L^{\frac{1}{1-s}}}(B_t)}=\lim\limits_{\epsilon\rightarrow 0^+}\|u^{-}_{\epsilon}\|_{L^{\frac{1}{1-s}}(B_t)}=0.
\end{eqnarray*}
Since $t$ is arbitrary in $(0,1)$, letting $t\rightarrow 1^-$, one has
\begin{eqnarray*}
\|u^-\|_{{L^{\frac{1}{1-s}}}(B_1)}=0.
\end{eqnarray*}
This implies that $u(x)\geq 0$ in $B_1$ and completes the proof of the theorem.
\end{proof}

\subsection{Maximum principles for fractional Laplacian with zero order term}
In this subsection, we prove the maximum principle and a refined version of the strong maximum principle for fractional Laplacian with zero order term(fractional Schr\"{o}dinger operator).

We first give the proof of Theorem \ref{thm7}.\\
\\
{\bf Proof of Theorem \ref{thm7}.}
\quad ~ It follows from Lemma \ref{lemma1} and (\ref{7}) that
\begin{eqnarray}\label{6.1}
 \left\{\begin{array}{l} (-\Delta)^s u^{-}(x)+c(x)u^{-}(x)\leq 0 \ \ \ \mbox{in} \ \mathcal{D}'(B_1)   \\
u^{-}(x)=0 \ \ \ \mbox{in}\ B_1^c.
\end{array}
\right.
 \end{eqnarray}

Let
\begin{eqnarray*}
v(x)=\left\{\begin{array}{l}\displaystyle{\int_{B_1}G_1(x,y)(-c(y)u^{-}(y))dy} \ \ \ \mbox{in} \ B_1  \\
0 \ \ \ \mbox{in}\ B_1^c.
\end{array}
\right.
\end{eqnarray*}
Note that
\begin{eqnarray*}
0<G_1(x,y)<\Phi(x-y) \ \ \text{for any }x,y\in B_1, x\neq y.
\end{eqnarray*}

Thus, for any $p>\frac{n}{n-2s}$, it follows from Hardy-Littlewood-Sobolev and H\"older inequalities that
\begin{eqnarray}\label{est}
\|v\|_{L^p(B_1)}&\leq& \|\int_{B_1}|\Phi(x-y)c^-(y)u^{-}(y)|dy\|_{L^p(B_1)}\\ \nonumber
                       &\leq& C(n,s,p)\|c^-u^{-}\|_{L^{\frac{np}{n+2sp}}(B_1)}\\ \nonumber
                       &\leq&  C(n,s,p)\|c^-\|_{L^{\frac{n}{2s}}(B_1)}\|u^{-}\|_{L^p(B_1)}.\nonumber
\end{eqnarray}
Clearly, $v(x)\in \mathcal{L}_{2s}\cap L^{\frac{1}{1-s}}(B_1)$. On the other hand, $v(x)$ satisfies the following equation
\begin{eqnarray}\label{6.2}
 \left\{\begin{array}{l} (-\Delta)^s v(x)+c(x)u^{-}(x)= 0 \ \ \ \mbox{in} \ \mathcal{D}'(B_1)   \\
v(x)=0 \ \ \ \mbox{in}\ B_1^c.
\end{array}
\right.
 \end{eqnarray}

Combining (\ref{6.1}) and (\ref{6.2}) yields that
\begin{eqnarray}\label{6.3}
 \left\{\begin{array}{l} (-\Delta)^s (v(x)-u^{-}(x))\geq 0 \ \ \ \mbox{in} \ \mathcal{D}'(B_1)   \\
v(x)-u^{-}(x)=0 \ \ \ \mbox{in}\ B_1^c.
\end{array}
\right.
 \end{eqnarray}
It follows from Theorem \ref{thm6} that $v(x)\geq u^{-}(x)$ in $B_1$.

Thus, from(\ref{est}), one can derive that
\begin{eqnarray*}
\|u^{-}\|_{L^{\frac{1}{1-s}}(B_1)}\leq  C(n,s)\|c^-\|_{L^{\frac{n}{2s}}(B_1)}\|u^{-}\|_{L^{\frac{1}{1-s}}(B_1)}.
\end{eqnarray*}

If
 $$\|c^-\|_{L^{\frac{n}{2s}}(B_1)}< \frac{1}{C(n,s)},$$
then
\begin{eqnarray*}
\|u^{-}\|_{L^{\frac{1}{1-s}}{(B_1)}}=0.
\end{eqnarray*}
This implies that $u^{-}(x)=0$ in $B_1$. Thus $u(x)\geq 0$ in $B_1$.

Hence the proof of the theorem is completed.
\hfill\qedsymbol
\\

In the following, we present the proof of the refined version of the strong maximum principle for fractional Laplacian with zero order term if $p>\frac{n}{2s}$. \\
\\
{\bf Proof of Theorem \ref{thms}.}
\quad ~ It follows from Theorem \ref{thm7} that
\begin{eqnarray*}
u(x)\geq 0 \ \ \mbox{in} \ B_1.
\end{eqnarray*}
Then (\ref{09}) can be written as
\begin{eqnarray}\label{3.9.1}
 \left\{\begin{array}{l} (-\Delta)^s u(x)+c^{+}(x)u(x)\geq c^{-}(x)u(x)\geq 0 \ \ \ \mbox{in} \ \mathcal{D}'(B_1)   \\
u(x) \geq m>0 \ \ \ \mbox{in}\ B_2\setminus B_1\\
u(x)\geq 0 \ \ \ \mbox{in} \ B_2^c.
\end{array}
\right.
 \end{eqnarray}

Let
\begin{eqnarray*}
 v(x)=\left\{\begin{array}{l} \int_{B_1}G_1(x,y)c^{+}(y)dy \ \ \ \mbox{in} \ B_1  \\
 0 \ \ \ \mbox{in}\  B_1^c.
\end{array}
\right.
 \end{eqnarray*}

By direct calculations and classical elliptic estimates, one has
\begin{eqnarray}\label{3.9.2}
 \left\{\begin{array}{l} (-\Delta)^s v(x)=c^{+}(x) \ \ \ \mbox{in} \ \mathcal{D}'(B_1)   \\
v(x) = 0 \ \ \ \mbox{in}\  B_1^c,
\end{array}
\right.
 \end{eqnarray}
and
\begin{eqnarray*}
\|{v}\|_{L^{\infty}(B_1)}
\leq C(n,s,p)\|{c^{+}}\|_{L^p(B_1)} \ \ \text{   for }p>\frac{n}{2s}.
\end{eqnarray*}

Let
 \begin{eqnarray}
 h(x)=\left\{\begin{array}{l} \int_{B_2\setminus B_1}P_1(x,y)mdy \ \ \ \mbox{in} \ B_1  \\
 m \ \ \ \mbox{in}\  B_2\setminus B_1\\
 0  \ \ \ \mbox{in} \ B_2^c.
\end{array}
\right.
 \end{eqnarray}
Then $h(x)$ satisfies the following problem
 \begin{eqnarray}\label{3.9.4}
 \left\{\begin{array}{l} (-\Delta)^s h(x)=0 \ \ \ \mbox{in} \ \mathcal{D}'(B_1)   \\
h(x) = m \ \ \ \mbox{in}\  B_2\setminus B_1\\
h(x) = 0  \ \ \ \mbox{in} \ B_2^c.
\end{array}
\right.
 \end{eqnarray}

 Moreover, there exists a positive constant $\beta<1$ depending on $n$, $s$ such that
 \begin{eqnarray}\label{3.9.5}
 \beta m \leq h(x) \leq m \ \ \mbox{in} \ B_1.
 \end{eqnarray}

 Define $w(x)=u(x)-h(x)+mv(x)$. It follows from (\ref{3.9.1})-(\ref{3.9.5}) that $w(x)$ satisfies the following equation
\begin{eqnarray}\label{3.3}
 \left\{\begin{array}{l} (-\Delta)^s w(x)+c^{+}(x)w(x)\geq 0 \ \ \ \mbox{in} \ \mathcal{D}'(B_1)   \\
w(x) \geq 0 \ \ \ \mbox{in}\  B_1^c.
\end{array}
\right.
 \end{eqnarray}
Thus, $w(x)\geq 0$ in $B_1$ by Theorem \ref{thm7}.

Therefore, one has, for $x\in B_1$,
\begin{eqnarray*}
u(x)&\geq& h(x)-mv(x)\\
&\geq&  \beta m- m\|v\|_{L^{\infty}(B_1)}\\
&\geq &  m(\beta- C(n,s,p)\|c^{+}\|_{L^p(B_1)}).
\end{eqnarray*}

Let
$$k(n,s,p)=\frac{\beta}{2C(n,s,p)}.$$
If $\|{c^{+}}\|_{L^p(B_1)}\leq k(n,s,p)$, then one has
$$u(x)\geq \frac{\beta}{2}m.$$

Hence the proof of the theorem is completed.
\hfill\qedsymbol

\subsection{Maximum principles for fractional Laplacian with first order term}
In this subsection, we prove the maximum principle and strong maximum principle for fractional Laplacian with first order term.

We first give the proof of Theorem \ref{thm8}.\\
\\
{\bf Proof of Theorem \ref{thm8}.}~~
\quad It follows from Lemma \ref{lemma1} and (\ref{8}) that
\begin{eqnarray}\label{8.1}
 \left\{\begin{array}{l} (-\Delta)^s u^{-}(x)+\vec{b}(x)\cdot \nabla u^{-}(x)\leq 0 \ \ \ \mbox{in} \ \mathcal{D}'(B_1)  \\
u^{-}(x)=0 \ \ \ \mbox{in}\ B_1^c.
\end{array}
\right.
 \end{eqnarray}

Let
\begin{eqnarray*}
v(x)=\left\{\begin{array}{l}\displaystyle{\int_{B_1}div_y(G_1(x,y)\vec{b}(y))u^{-}(y)dy} \ \ \ \mbox{in} \ B_1  \\
0 \ \ \ \mbox{in}\ B_1^c.
\end{array}
\right.
\end{eqnarray*}
Thus, for $x\in B_1$, one has
\begin{eqnarray*}
v(x)&=&\int_{B_1}\nabla_y G_1(x,y)\cdot \vec{b}(y) u^-(y)dy+\int_{B_1}(div_y\vec{b}(y))G_1(x,y)u^-(y)dy\\
    &:=&T_1(x)+T_2(x).
\end{eqnarray*}

In the following, we estimate $T_1(x)$ and $T_2(x)$, respectively.

A simple calculation yields that
\begin{eqnarray*}
\nabla_yG_1(x,y)\!=\!\kappa(n,s)\!\left(\frac{(2s-n)(y-x)}{|x-y|^{n-2s+2}}\!\int_{0}^{\sigma_1(x,y)}\!\!\frac{t^{s-1}}{(t+1)^{\frac{n}{2}}}dt\!+\!\frac{1}{|x-y|^{n-2s}}\frac{\sigma_1^{s-1}(x,y)}{(\sigma_1(x,y)+1)^{\frac{n}{2}}}\nabla_y\sigma_1(x,y)\right).
\end{eqnarray*}
Note that
\begin{eqnarray*}
\nabla_y\sigma_1(x,y)&=&(1-|x|^2)\left(\frac{-2y}{|x-y|^2}+\frac{2(x-y)(1-|y|^2)}{|x-y|^4}\right)\\
              &=&\left(\frac{-2y}{1-|y|^2}+\frac{2(x-y)}{|x-y|^2}\right)\sigma_1(x,y).
\end{eqnarray*}

Thus, one has
\begin{eqnarray*}
&&|\nabla_yG_1(x,y)|\\
&&\!\!\!\!\!\!\!\!\leq C(n,s)\left(\frac{1}{|x-y|^{n-2s+1}}\int_{0}^{\sigma_1(x,y)}\frac{t^{s-1}}{(t+1)^{\frac{n}{2}}}dt+\frac{1}{|x-y|^{n-2s}}\left(\frac{1}{1-|y|}+\frac{1}{|x-y|}\right)\frac{\sigma_1^s(x,y)}{(\sigma_1(x,y)+1)^{\frac{n}{2}}}\right)\\
&&\!\!\!\!\!\!\!\!\leq C_1(n,s)\left(\frac{1}{|x-y|^{n-2s+1}}+\frac{1}{|x-y|^{n-2s}}\frac{1}{1-|y|}\right).
\end{eqnarray*}

Then
\begin{eqnarray}\label{T1}
|T_1(x)|\leq C_1(n,s)\int_{B_1}\left(\frac{1}{|x-y|^{n-2s+1}}+\frac{1}{|x-y|^{n-2s}}\frac{1}{1-|y|}\right)|\vec{b}(y)|u^-(y)dy.
\end{eqnarray}

On the other hand, note that
\begin{eqnarray*}
0<G_1(x,y)<\Phi(x-y) \ \ \text{for any }x,y\in B_1, x\neq y.
\end{eqnarray*}
Thus, one has
\begin{eqnarray}\label{T2}
|T_2(x)|\leq C_2(n,s)\int_{B_1}\frac{|(div_y\vec{b}(y))u^-(y)|}{|x-y|^{n-2s}}dy.
\end{eqnarray}

Therefore, for $p>\frac{n}{n-2s}$, combining (\ref{T1}) and (\ref{T2}), together with {\color{red}} Hardy-Littlewood-Sobolev and H\"older inequalities yields
\begin{eqnarray}\label{v}
\|v\|_{L^p(B_1)}\!\!\!\!&\leq & C_3(n,s)\|\int_{B_1}\left\{\frac{|\vec{b}(y)|u^-(y)}{|x-y|^{n-2s+1}}+\frac{u^-(y)}{|x-y|^{n-2s}}\left(\frac{|\vec{b}(y)|}{d(y)}+|div_y\vec{b}(y)|\right)\right\}dy\|_{L^p(B_1)}\\ \nonumber
                &\leq & C(n,s,p)\left(\|\vec{b}u^{-}\|_{L^{\frac{np}{n+(2s-1)p}}(B_1)}+\|(div\vec{b})u^{-}\|_{L^{\frac{np}{n+2sp}}(B_1)}+\|\frac{\vec{b}}{d}u^{-}\|_{L^{\frac{np}{n+2sp}}(B_1)}\right)\\ \nonumber
                &\leq &  C_1(n,s,p)\left(\|\vec{b}\|_{L^{\frac{n}{2s-1}}(B_1)}+\|div\vec{b}\|_{L^{\frac{n}{2s}}(B_1)}+\|\frac{\vec{b}}{d}\|_{L^{\frac{n}{2s}}(B_1)}\right)\|u^{-}\|_{L^p(B_1)}\\ \nonumber
                &\leq & C_2(n,s,p)\left(\|\vec{b}\|_{W^{1,\frac{n}{2s}}(B_1)}+\|\frac{\vec{b}}{d}\|_{L^{\frac{n}{2s}}(B_1)}\right)\|u^{-}\|_{L^p(B_1)}.\nonumber
\end{eqnarray}
 Clearly, $v(x)\in \mathcal{L}_{2s}\cap L^{\frac{1}{1-s}}(B_1)$. On the other hand, $v(x)$ satisfies the following equation
\begin{eqnarray}\label{8.2}
 \left\{\begin{array}{l} (-\Delta)^s v(x)+\vec{b}(x)\cdot \nabla u^{-}(x)= 0 \ \ \ \mbox{in} \ \mathcal{D}'(B_1)  \\
v(x)=0 \ \ \ \mbox{in}\ B_1^c.
\end{array}
\right.
 \end{eqnarray}

Combining (\ref{8.1}) and (\ref{8.2}) yields
\begin{eqnarray}\label{8.3}
 \left\{\begin{array}{l} (-\Delta)^s (v(x)-u^{-}(x))\geq 0 \ \ \ \mbox{in} \ \mathcal{D}'(B_1)  \\
v(x)-u^{-}(x)=0 \ \ \ \mbox{in}\ B_1^c.
\end{array}
\right.
 \end{eqnarray}
It follows from Theorem \ref{thm6} that $v(x)\geq u^{-}(x)$ in $B_1$.

Thus, from (\ref{v}), one can derive that
\begin{eqnarray*}
\|u^{-}\|_{L^{\frac{1}{1-s}}(B_1)}\leq  C_4(n,s)\left(\|\vec{b}\|_{W^{1,\frac{n}{2s}}(B_1)}+\|\frac{\vec{b}}{d}\|_{L^{\frac{n}{2s}}(B_1)}\right)\|u^{-}\|_{L^{\frac{1}{1-s}}(B_1)}.
\end{eqnarray*}

If
\begin{eqnarray*}
\|\vec{b}\|_{W^{1,\frac{n}{2s}}(B_1)}+\|\frac{\vec{b}}{d}\|_{L^{\frac{n}{2s}}(B_1)}<\frac{1}{C_4(n,s)},
\end{eqnarray*}
then
\begin{eqnarray*}
\|u^{-}\|_{L^{\frac{1}{1-s}}(B_1)}=0.
\end{eqnarray*}
Therefore, one has $u^{-}(x)=0$ in $B_1$. This implies that $u(x)\geq 0$ in $B_1$.
\hfill\qedsymbol
\\

Finally, we show the strong maximum principle for fractional Laplacian with first order term.\\
\\
{\bf Proof of Theorem \ref{thm9}.}~~\quad Let
\begin{eqnarray}
l(x)=\frac{\int_{B_2\setminus B_1}\frac{1}{|x-y|^{n+2s}}dy}{\int_{B_2\setminus B_1}\frac{1}{|x-y|^{n+2s}}dy+\int_{B_2^c}\frac{1}{|x-y|^{n+2s}}dy},\ \ \ \ x\in B_1.
\end{eqnarray}
By direct calculations, there exist constants $C_1(n,s)$ and $C_2(n,s)$ such that
\begin{eqnarray*}
\int_{B_2\setminus B_1}\frac{1}{|x-y|^{n+2s}}dy\geq C_1(n,s) \ \ \mbox{and} \ \ \int_{B_2^c}\frac{1}{|x-y|^{n+2s}}dy\leq C_2(n,s) \ \mbox{ for  } x\in B_1.
\end{eqnarray*}

Let $\gamma=\frac{C_1(n,s)}{C_1(n,s)+C_2(n,s)}$. Thus, one has
\begin{eqnarray}
l(x)\geq \gamma,\ \ \ x\in B_1.
\end{eqnarray}

Define
\begin{eqnarray}
 h(x)=\left\{\begin{array}{l}  \gamma m \ \ \ \mbox{in}\  B_1\\
 m\ \ \ \mbox{in}\  B_2\setminus B_1\\
 0  \ \ \ \mbox{in} \ B_2^c,
\end{array}
\right.
 \end{eqnarray}
and
 $$v(x)=u(x)-h(x).$$

One can derive that
\begin{eqnarray}\label{5.1}
 \left\{\begin{array}{l} (-\Delta)^s v(x)+\vec{b}(x)\cdot \nabla v(x)\geq 0 \ \ \ \ \mbox{in} \ \mathcal{D}'(B_1) \\
v(x) \geq 0 \ \ \ \mbox{in} \ B_1^c.
\end{array}
\right.
 \end{eqnarray}
 It follows from Theorem \ref{thm8} that $v(x) \geq 0$ in $B_1$. Thus,
\begin{eqnarray*}
u(x)\geq \gamma m \ \ \mbox{in}\ B_1.
\end{eqnarray*}

Hence the proof of the theorem is completed.

\section*{Acknowledgments} The research of Li was partially supported by NSFC grants 12031012 and 11831003. The authors would like to thank Professor Genggeng Huang and Chenkai Liu for their helpful discussions.

\end{document}